\documentclass[11pt]{amsart}
\usepackage{amssymb}
\usepackage{epsfig}
\usepackage{graphicx}

\usepackage{fullpage}
\usepackage[numbers,sort&compress]{natbib}
\usepackage{color}
\usepackage{setspace}

\def\endproof{\hfill$\square$\vspace{6pt}}

\def\alp{\alpha}
\def\eps{\varepsilon}
\newcommand{\R}{{\mathbb R}}
\begin{document}
\newtheorem*{proposition*}{Proposition~$3.7^{(')}$}
\newcommand{\myref}[1]{7$'$}
\newtheorem{theorem}{Theorem}
\newtheorem{definition}{Definition}
\newtheorem{corollary}{Corollary}
\newtheorem{proposition}{Proposition}
\newtheorem{lemma}{Lemma}
\newtheorem{remark}{Remark}
\newtheorem{example}{Example}
\newtheorem{assumption}{Assumption}
\newtheorem{conjecture}{Conjecture}
\newtheorem{problem}{Problem}

\title[]{Classification of solutions to equations involving Higher-order fractional laplacian}

\author{Zhuoran Du }
\address{School of Mathematics, Hunan University, Changsha 410082,
                 PRC  }
\email{duzr@hnu.edu.cn}

\author{Zhenping Feng }
\address{School of Mathematics, Hunan University, Changsha 410082,
                 PRC  }
\email{fengzp@hnu.edu.cn}

\author{Jiaqi Hu }
\address{School of Mathematics, Hunan University, Changsha 410082,
                 PRC  }
\email{13592528691@163.com}

\author{Yuan Li }
\address{School of Mathematics, Hunan University, Changsha 410082,
                 PRC  }
\email{liy93@hnu.edu.cn}

\thanks{}

\date{}

\maketitle

\begin{abstract}
In this paper, we are concerned with the following equation involving higher-order fractional Lapalacian
\begin{equation*}
\left\{\begin{aligned}
&(-\Delta)^{p+{\frac{\alpha}{2}}}u(x)=u_+^\gamma~~ \mbox{ in }\mathbb{R}^n,\\
&\int_{\mathbb{R}^n}u_+^\gamma dx<+\infty,
\end{aligned}\right.
\end{equation*}
where $p\geq 1$ is an integer, $0<\alp<2$,  $n> 2p+\alpha$ and $\gamma \in (1,\frac{n}{n-2p-\alp})$. We establish an integral representation formula for any nonconstant classical solution satisfying certain growth at infinity. From this we prove that these solutions  are  radially  symmetric about some point in $\R^n$ and monotone decreasing in the radial direction via method of moving planes in integral forms.

\end{abstract}


\section{Introduction}

In this paper, we study the following higher-order fractional equation
\begin{equation}\label{T1}
\left\{\begin{aligned}
&(-\Delta)^{p+\frac{\alp}{2}}u(x)=u_+^\gamma,\ \ \ \ x\in\mathbb{R}^n,
\\&\int_{\mathbb{R}^n}u_+^\gamma dx<+\infty,
\end{aligned}\right.
\end{equation}
where $\gamma \in (1,\frac{n}{n-2p-\alp})$, $n> 2p+\alpha$, $0<\alp<2$, $p\geq 1$ is an integer, $u_+=\max\{u,0\}$ and the higher-order fractional Laplacian is defined by $(-\Delta)^{p+\frac{\alp}{2}}:=(-\Delta)^{p}(-\Delta)^{\frac{\alp}{2}}$. It is well-known that the integral form of fractional Laplacian $(-\Delta)^{\frac{\alpha}{2}}$ $(0<\alpha<2)$ is
$$(-\Delta)^{\frac \alp 2}u(x)=C_{n,\alp}P.V.\int_{\R^n}\frac{u(x)-u(y)}{|x-y|^{n+\alp}}dy,$$
where the constant $C_{n, \alpha}=\frac{2^{\alpha}\Gamma(\frac{n+\alpha}{2})}{-\pi^{n/2}\Gamma(-\frac{\alpha}{2})}$.

Denote
$$\mathcal{L}_\alp(\R^n):=\Big\{u:\R^n\to\R \big|\int_{\R^N}\frac{|u(x)|}{1+|x|^{n+\alp}}dx<+\infty\Big\}.$$
It is known that $(-\Delta)^{\frac \alp 2}u$ is well defined for $u\in C_{loc}^{[\alp],\{\alp\}+\eps}\cap \mathcal{L}_\alp(\R^n)$, where $[\alp]$ denotes the integer part of $\alp$, $\{\alp\}:=\alp-[\alp]$ and  any $\eps>0$. In order to guarantee that $(-\Delta)^{\frac \alp 2} u\in C^{2p}(\R^n)$, we have to assume $u\in C_{loc}^{2p+[\alp],\{\alp\}+\eps}\cap \mathcal{L}_\alp(\R^n)$ (see \cite{16}, \cite{Ch}), and hence $u$ is a classical solution of the equation in (\ref{T1}) in the sense that $(-\Delta)^{p+\frac \alp 2}u$ is point-wise well-defined in the whole $\mathbb{R}^n$.

For decades, many researchers are interested in the classification of solutions to semi-linear elliptic equations. We now recall some relevant results. The following conformally invariant equations
\begin{equation}\label{1.3}
(-\Delta)^{\frac \beta 2}u=u^{\frac{n+\beta}{n-\beta}} ~\mbox{in} ~\mathbb{R}^n,~0<\beta<n
\end{equation}
have been extensively studied (see \cite{2,S,9,10,DQ,Dai,13,15,17} and the reference therein). In case $\beta=2$, (\ref{1.3}) becomes the well-known Yamabe equation. Cafferelli, Gidas and Spruck \cite{2} (see also  \cite{7}) classified all positive solutions of Yamabe equation. For general $0<\beta<n$, Chen, Li and Ou in \cite{10} (see also \cite{17}) classified all positive $L_{loc}^{\frac{2n}{n-\beta}}$ solutions to the equivalent integral equation of (\ref{1.3}) by the method of moving planes in integral forms. For $0<\beta< 2$, Chen, Li and Li in \cite{9} developed a direct method of moving planes for fractional Laplace equation and complete the classification of all nonnegative solutions. When $\beta = 3$ is an odd integer, Dai and Qin in \cite{Dai} derived the classification of nonnegative classical solutions to (\ref{1.3}) with the assumption $\int_{\R^n}\frac{u^\frac{n+3}{n-3}}{|x|^{n-3}}dx<+\infty$. Recently, Cao, Dai and Qin in \cite{3} extended the result of \cite{Dai} to general case $0<\beta<n$.
Their classification results completely improved the results in \cite{Dai} without assumption on integrability. Precisely they proved the super poly-harmonic properties for nonnegative solutions by making full use of the the Poisson representation formula for $(-\Delta)^{\frac \alp 2}(0<\alpha<2)$ and developing some new integral estimates on the outer-spherical average $\int_R^{+\infty}\frac{R^\alp}{r(r^2-R^2)^{\frac \alp 2}}dr$ and iteration techniques.
Based on super poly-harmonic properties, they established integral representation formula, Liouville's Theorem and classification of classical solutions to higher-order equations involving fractional Laplacian.

 For (\ref{T1}), classification of nonconstant solutions have been established for some particular $p, \alpha$ and $\gamma$. For $p=1, \alpha=0$  and $\gamma=\frac{n}{n-2}$, Wang and Ye \cite{Wang} classified  all nonconstant solutions of (\ref{T1}).
Suzuki and Takahashi \cite{Suzu} extended the results of \cite{Wang} from the exponent $\frac{n}{n-2}$  to  more general exponent. Precisely, they considered the problem
\begin{equation}\label{0a1}
-\Delta v=v^{\gamma}_+  ~~ \mbox{ in }{\mathbb{R}^n},~~n>2,~~\int_{\mathbb{R}^n}v_+^{\frac{n(\gamma-1)}{2}}dx<+\infty,
\end{equation}
where $\gamma\in(1,\frac{n+2}{n-2})$. For the case $p=2, \alpha=0$  and $\gamma\in (1,\frac{n}{n-4}]$, Chammakhi, Harrabi and Selmi \cite{4} completed the classification of all sign-changing solutions  of (\ref{T1}). The first, second and fourth authors \cite{Du} extended the result of \cite{4} to the corresponding polyharmonic equation, namely for the case of any integer $p$ larger than $1$. Particularly for the case $p=2$, they obtain the same results for more general exponent $\gamma\in(1,\frac{n+4}{n-4})$. They also considered nonconstant solutions of the following fractional equation (namely $p=0$ corresponds to the equation in (\ref{T1}))
\begin{equation}\label{p7}
(-\Delta)^{\frac{\alpha}{2}}v=v^\gamma_+  ~~ \mbox{ in }{\mathbb{R}^n},~~\alpha\in(0,2),~~n\geq 2,~~
\int_{\mathbb{R}^n}v_+^{\frac{n(\gamma-1)}{\alpha}}dx<+\infty,~~~~~
\end{equation}
where $\gamma \in (1, \frac{n+\alpha}{n-\alpha})$, and obtained the corresponding classification results.

To classify nonconstant solutions of (\ref{T1}) for general $0<\alp<2$,  $p\geq 1$ and $\gamma \in (1,\frac{n}{n-2p-\alp})$, we need to establish the following super poly-harmonic properties.

\begin{theorem}\label{Th1}
Assume $\gamma \in(1,\frac{n}{n-2p-\alp})$, $n>2p+\alp$, $0<\alp<2$ and $p\geq 1$ is an integer. Suppose $u$ is a nonconstant classical solution to (\ref{T1}) satisfying $u(x)=o(|x|^{\alp+\frac 4 \alp})$. Then for every $i=0,1,\ldots,p-1$,
$$(-\Delta)^{i+\frac{\alp}{2}}u(x)\geq 0,\ \ \ \ \forall~ x\in \R^n.$$
\end{theorem}
As well-known, establishing equivalent integral equation of (1.1) is beneficial to classifying solutions to higher order equations. As a consequence of super poly-harmonic properties for equation (\ref{T1}), we can  derive the equivalence between the differential equation (\ref{T1}) and the following integral equation
\begin{equation}\label{1.2}
u(x)=\int_{\R^n}\frac{R_{2p+\alp,n}}{|x-y|^{n-2p-\alp}}u_+^\gamma (y)dy+C_0,
\end{equation}
where $C_0<0$ and the Riesz potential's constant $R_{\mu,n}=\frac{\Gamma(\frac{n-\mu}{2})}{\pi^{\frac n 2}2^\mu\Gamma(\frac\mu 2)}$ for $0<\mu<n$ (see \cite{MS}).
\begin{theorem}\label{T2}
Under the  assumptions of Theorem \ref{Th1}, suppose u is a nonconstant classical solution to (\ref{T1}) satisfying $u(x)=o(|x|^{\alp+\frac 4 \alp})$. Then u is a solution to integral equation (\ref{1.2}), and vice verse. Moreover, the support of $u_+$ is compact.
\end{theorem}

Apparently, for (\ref{T1}), any non-positive constant is its solution and any positive constant is not its solution.  Hence we only consider nonconstant solutions. Actually, there is no positive solutions to (\ref{T1}) from the nonexistence of positive  solutions in the subcritical exponent case of  Theorem 1.9 in \cite{3}. All negative nonconstant solutions of (\ref{T1}) can be ruled out by the Liouville theorem for fractional poly-harmonic functions (see Theorem 1.3 in \cite{3}). Consequently any nonconstant solutions of (\ref{T1}) are sign-changing solutions. From this observation  and Theorem \ref{T2}, the following classification results can be obtained by applying the method of moving planes in integral forms.
\begin{theorem}\label{T3}
Under the  assumptions of Theorem \ref{Th1}, suppose $u$ is a sign-changing classical solution to equation (\ref{T1}) satisfying $u(x)=o(|x|^{\alp+\frac 4 \alp})$ at infinity, then $u$ is symmetric about some point $x_0\in \R^n$ and $\frac{\partial u}{\partial r}<0$, where $r=|x-x_0|$.
\end{theorem}

Throughout this paper, we will use $C$ to denote a positive constant which may change  from line to line and even within the same line.

\section{Preliminaries}
The purpose of this section is to introduce several useful properties which will be crucial to the forthcoming sections.
\begin{proposition}\label{o1}\cite{B}
Let $R>0$, $x_0\in \R^n$, $h\in C^{\alp+\eps}(B_R(x_0))\cap C(\overline{B}_R(x_0))$, $0<\alp<2$ and let
$$u(x):=\left\{\begin{aligned}
&\int_{B_R(x_0)}G_R^\alp (x,y)h(y)dy,&\ \ \ \ &x\in B_R(x_0),
\\&0,&\ \ \ \ &x\in \R^n\setminus B_R(x_0).
\end{aligned}\right.
$$
Then $u$ is the unique point-wise continuous solution of the following problem
$$\left\{\begin{aligned}
&(-\Delta)^{\frac \alp 2}u(x)=h(x),&\ \ \ \ &x\in B_R(x_0),
\\&u(x)=0,&\ \ \ \ &x\in \R^n\setminus B_R(x_0).
\end{aligned}\right.
$$
The Green function is defined by
$$G_R^\alp(x,y):=\frac{C_{n,\alp}}{|x-y|^{n-\alp}}\int_0^{\frac{t_R}{s_R}}\frac{b^{\frac{\alp}{2}-1}}{(1+b)^{\frac{n}{2}}}db,\ \ \ \ \  x,y\in B_R(x_0)$$
with $s_R=\frac{|x-y|^2}{R^2}$, $t_R=\Big(1-\frac{|x-x_0|^2}{R^2}\Big)\Big(1-\frac{|y-x_0|^2}{R^2}\Big)$, and $G_R^\alp(x,y)=0$ if $x$ or $y\in \R^n\setminus B_R(x_0)$,
where
$$C_{n,\alp}=\frac{\Gamma(\frac n 2)}{2^{\alp}\pi^{\frac n 2}\Gamma^2(\frac \alp 2)}\ \mbox{if} \ \ n\neq \alp,\ \ \ \ C_{1,1}=\frac 1\pi\  \mbox{if} \ n=\alp=1.$$
\end{proposition}

\begin{proposition}\label{o2}\cite{B}
Let $R>0$, $x_0\in \R^n$, $g\in\mathcal{L}_\alp(\R^n)\cap C(\R^n)$, $0<\alp<2$ and let
$$u_g(x):=\left\{\begin{aligned}
&\int_{\R^n\setminus B_R(x_0)}P_R^\alp(x,y)g(y)dy,&\ \ \ \ &x\in B_R(x_0),
\\&g(x),&\ \ \ \ &x\in \R^n\setminus B_R(x_0).
\end{aligned}\right.
$$
Then $u_g$ is the unique point-wise continuous solution of the following problem
$$\left\{\begin{aligned}
&(-\Delta)^{\frac \alp 2}u(x)=0,&\ \ \ \ &x\in B_R(x_0),
\\&u(x)=g(x),&\ \ \ \ &x\in \R^n\setminus B_R(x_0).
\end{aligned}\right.
$$
The Poisson kernel $P_R^\alp$ is defined by
$$P_R^\alp(x,y):=\frac{\Gamma(\frac{n}{2})}{\pi^{\frac{n}{2}+1}}\sin\frac{\pi\alp}{2}\Big(\frac{R^2-|x-x_0|^2}{|y-x_0|^2-R^2}\Big)^{\frac\alp 2}\frac 1{|x-y|^n}.$$
\end{proposition}

\begin{proposition}\label{o3}\cite{9}
Let $\Omega$ be a bounded domain in $\R^n$. Assume that $u\in\mathcal {L}_\alp\cap C_{loc}^{1,1}(\Omega)$ and is lower semi-continuous on $\bar{\Omega}$. If
$$\left\{\begin{aligned}
&(-\Delta)^{\frac \alp 2}u(x)\geq0,&\ \ \ \ &x\in \Omega,
\\&u(x)\geq 0,&\ \ \ \ &x\in \R^n\setminus \Omega,
\end{aligned}\right.
$$
then $u(x)\geq 0$ in $\Omega$. Moreover, if $u=0$ at some point in $\Omega$, then $u=0$ a.e. in $\R^n$. These conclusions also holds for unbounded domain $\Omega$, if we further assume that $\mathop{\liminf}\limits_{|x|\to +\infty }u(x)\geq 0$.
\end{proposition}

\section{Proof of Theorem \ref{Th1}}

In this section we will complete the proof of Theorem \ref{Th1}.

{\bf Proof} We borrow the idea in \cite{3} to prove this theorem. Denote $v_i:=(-\Delta)^{i+\frac{\alp}{2}}u(x)$ for $i=0,1,\ldots,p-1$. It follows from (\ref{T1}) that
\begin{equation}\label{2.1}
\left\{\begin{aligned}
&(-\Delta)^{\frac{\alp}{2}}u(x)=v_0,&\ \ \ \ &x\in\mathbb{R}^n,
\\&-\Delta v_0=v_1,&\ \ \ \ &x \in\mathbb{R}^n,
\\&\cdots\cdots
\\&-\Delta v_{p-1}=u_+^\gamma,&\ \ \ \ & x\in\mathbb{R}^n.
\end{aligned}\right.
\end{equation}
Assume that Theorem \ref{Th1} is not true, then there must exists a largest integer $0\leq k\leq p-1$ and a point $x_0 \in\R^n$ such that
$$v_k(x_0)=(-\Delta)^{k+\frac{\alp}{2}}u(x_0)<0.$$
For any $r>0$, define
$$\bar{g}(r,x_0):=\frac{1}{|\partial B_r(x_0)|}\int_{\partial B_r(x_0)}g(x)d\sigma,$$
where $|\partial B_r(x_0)|$ denotes the area of the sphere $\partial B_r(x_0)$. For simplicity, we denote $\bar{g}(r,x_0)$ as $\bar{g}(r)$.

First, we will illustrate that $0\leq k\leq p-1$ is even by contradiction. Assume $k$ is odd. From the well-known property $\overline{\Delta u}(r)=\frac{1}{r^{n-1}}(r^{n-1}\bar{u}'(r))'$ and (\ref{2.1}), we have
$$\bar{v}_k(r)\leq\bar{v}_k(0):=-c_0<0, \ \ \ \ \forall r>0.$$
Simple calculation shows that
$$\bar{v}_{k-1}(r) \geq\bar{v}_{k-1}(0)+\frac{c_0}{2n}r^2, \ \ \ \ \forall r>0,$$
and
$$\bar{v}_{k-2}(r)\leq\bar{v}_{k-2}(0)-\frac{\bar{v}_{k-1}(0)}{2n}r^2-\frac{c_0}{8n(n+2)}r^4,
 \ \ \ \ \forall r>0.$$
Repeating the above argument, we derive
\begin{equation}\label{2.2}
\bar{v}_{0}(r)\geq\bar{v}_{0}(0)+c_1r^2+c_2r^4+\cdots+c_kr^{2k},\ \ \ \ \forall r>0,
\end{equation}
where $c_k>0$. From (\ref{2.2}), we deduce that that there exists a $r_0$ large enough such that
\begin{equation}\label{2.3}
\bar{v}_{0}(r)\geq \frac{1}{2}c_kr^{2k}\ \ \ \ \forall r>r_0.
\end{equation}

The first equation in (\ref{2.1}), Proposition \ref{o1} and Proposition \ref{o2} imply that for arbitrary $R>0$,
$$u(x)=\int_{B_R(x_0)}G_R^\alp(x,y)v_0(y)dy+\int_{|y-x_0|>R}P_R^\alp(x,y) u(y)dy,\ \ \ \ \forall x \in B_R(x_0).$$
Therefore, we have
\begin{align}\label{2.4}
+\infty >u(x_0)&=\int_{B_R(x_0)}\frac{C_{n,\alp}}{|x_0-y|^{n-\alp}}\Big(\int_0^{\frac{R^2}{|y-x_0|^2}-1}
\frac{b^{\frac{\alp}{2}-1}}{(1+b)^{\frac{n}{2}}}db\Big)v_0(y)dy
\notag\\& +C'_{n,\alp}\int_{|y-x_0|>R}\frac{R^\alp}{(|y-x_0|^2-R^2)^{\frac\alp 2}} \frac {u(y)}{|y-x_0|^n}dy
\\&=:\mathrm{I}+\mathrm{II}.\nonumber
\end{align}
Observe that $0<r\leq\frac R2$ implies $3\leq\frac{R^2}{r^2}-1<+\infty$. Thus
\begin{equation}\label{2.5}
\int_0^{3}\frac{b^{\frac{\alp}{2}-1}}{(1+b)^{\frac{n}{2}}}db
\leq\int_0^{\frac{R^2}{r^2}-1}\frac{b^{\frac{\alp}{2}-1}}{(1+b)^{\frac{n}{2}}}db
\leq\int_0^{+\infty}\frac{b^{\frac{\alp}{2}-1}}{(1+b)^{\frac{n}{2}}}db.
\end{equation}
From (\ref{2.3})-(\ref{2.5}), we conclude that for any $R>2r_0$ there exist $c_i>0$ $(i=1,...,6)$ independent of $R$ such that
\begin{align}\label{2.6}
\mathrm{I}&=C_{n,\alp}|\partial B_1|\int_0^{R}r^{\alp-1}\Big(\int_0^{\frac{R^2}{r^2}-1}
\frac{b^{\frac{\alp}{2}-1}}{(1+b)^{\frac{n}{2}}}db\Big)\bar{v}_0(y)dy
\notag\\&\geq c_1\int_{r_0}^{\frac R 2}r^{\alp-1}\bar{v}_0(r)dr-c_2\int_0^{r_0}r^{\alp-1}|\bar{v}_0(r)|dr
\notag\\&\geq c_3\int_{r_0}^{\frac R 2}r^{2k+\alp-1}dr-c_4
\notag\\&\geq c_5R^{2k+\alp}-c_6.
\end{align}
Due to $u(x)=\mbox{o}(|x|^{\alp+\frac 4 \alp})$ at infinity, we derive that  for sufficiently large $R$
\begin{align}\label{E1}
C_{n,\alp}'&\int_{R<|y-x_0|<R+R^{1-\frac{4}{\alp}}}\frac{R^\alp}{(|y-x_0|^2-R^2)^{\frac\alp 2}} \frac {u(y)}{|y-x_0|^n}dy
\notag\\&=CR^\alp\int_R^{R+R^{1-\frac{4}{\alp}}}\frac{\bar{u}(r)}{r(r^2-R^2)^{\frac\alp 2}}dr
\notag\\&\leq \frac{2-\alpha}{2^{\frac{\alpha}{2}+\frac{4}{\alpha}+2}}c_5R^\alp\int_R^{R+R^{1-\frac{4}{\alp}}}\frac{r^{\alp+\frac 4 \alp}}{r(r^2-R^2)^{\frac\alp 2}}  dr
\notag\\&\leq \frac{2-\alpha}{8}c_5R^{\frac{3}{2}\alp+\frac{4}{\alp}-1}\int_R^{R+R^{1-\frac{4}{\alp}}}(r-R)^{-\frac{\alp}{2}}dr
\notag\\&\leq \frac{c_5}{4}R^{2+\alp}.
\end{align}
Note that if $|y-x_0|\geq R+R^{1-\frac 4 \alp}$, one has
$$\frac{R^2}{|y-x_0|^2}\leq \frac{R^2}{(R+R^{1-\frac{4}{\alp}})^2}\leq\frac{R^2}{R^2+2R^{2-\frac{4}{\alp}}},$$
which gives that
$$|y-x_0|^2-R^2\geq\frac{2R^{2-\frac 4 \alp}}{R^2+2R^{2-\frac 4 \alp}}|y-x_0|^2. $$
Obviously  for large $R$, $|y-x_0|\geq\frac{|y|}{2}$ holds if $|y-x_0|\geq R+R^{1-\frac 4 \alp}$. Therefore for $R$ large enough, we have
\begin{align}\label{E2}
C_{n,\alp}'&\int_{\R^n\setminus B_{R+R^{1-\frac 4 \alp}}(x_0)}\frac{R^\alp}{(|y-x_0|^2-R^2)^{\frac\alp 2}} \frac {u(y)}{|y-x_0|^n}dy
\notag\\&\leq CR^{\alp}\int_{\R^n\setminus B_{R+R^{1-\frac 4 \alp}}(x_0)}\frac{(R^2+2R^{2-\frac 4 \alp})^{\frac \alp 2}|u(y)|}{(R^{2-\frac 4 \alp })^{\frac \alp 2}|y|^{n+\alp}}dy
\notag\\&\leq CR^{2\alp-(\alp-2)}\int_{\R^n\setminus B_{R+R^{1-\frac 4 \alp}}(x_0)}\frac {|u(y)|}{|y|^{n+\alp}}dy
\notag\\&\leq \frac{c_5}{4} R^{\alp+2},
\end{align}
where the last inequality holds due to $u \in\mathcal{L}_\alp(\R^n)$. Combining (\ref{E1}) with (\ref{E2}), we conclude
\begin{equation}\label{2.8}
\mathrm{II}\leq \frac{c_5}{2} R^{\alp+2}.
\end{equation}
From (\ref{2.4}), (\ref{2.6}) and (\ref{2.8}), we derive that for $R$ large enough,
\begin{equation}\label{2.9}
+\infty>u(x_0)\geq c_5R^{2k+\alp}-c_6-\frac{c_5}{2} R^{2+\alp}.
\end{equation}
 Letting $R\to+\infty$ in (\ref{2.9}), we get immediately a contradiction. Thus, $k$ must be even.

Next, we will show that $k=0$. Otherwise, suppose that $2\leq k \leq p-1$ is even, using the same arguments as in deriving (\ref{2.2}), we deduce
\begin{equation}\label{2.10}
\bar{v}_{0}(r)\leq\bar{v}_{0}(0)-c_1r^2-c_2r^4-\cdots-c_kr^{2k},\ \ \ \ \forall r>0,
\end{equation}
where $c_k>0$. From (\ref{2.10}), we infer that there exists $r_1$ large enough such that
\begin{equation}\label{2.11}
\bar{v}_{0}(r)\leq -\frac{1}{2}c_kr^{2k}\ \ \ \ \forall r>r_1.
\end{equation}
Note that if $\frac{R} 2 <r<R$, then $0<\frac{R^2}{r^2}-1<3$, and hence
\begin{equation}\label{2.12}
\int_0^{\frac{R^2}{r^2}-1}\frac{b^{\frac{\alp}{2}-1}}{(1+b)^{\frac{n}{2}}}db\geq
\int_0^{\frac{R^2}{r^2}-1}\frac{b^{\frac{\alp}{2}-1}}{2^n}db\geq C \Big(\frac{R^2}{r^2}-1\Big)^{\frac \alp 2}.
\end{equation}
From (\ref{2.4}), (\ref{2.5}), (\ref{2.11}) and (\ref{2.12}), we get for any $R>2r_1$,
\[
\begin{split}
\int_R^{+\infty}&\frac{|\partial B_1|R^\alp \bar{u}(r)}{r(r^2-R^2)^{\frac\alp 2}}dr=-C_{n,\alp}\int_0^R r^{\alp-1}\Big(\int_0^{\frac{R^2}{r^2}-1}\frac{b^{\frac{\alp}{2}-1}}{(1+b)^{\frac{n}{2}}}db\Big)\bar{v}_0(r)dr+u(x_0)
\\&\geq C\int_{r_1}^{\frac R 2}r^{2k+\alp-1}dr-\widetilde{C}\int_0^{r_1}r^{\alp-1}|\bar{v}_0(r)|dr+C\int_{\frac R 2}^R r^{2k+\alp-1}\Big(\frac{R^2}{r^2}-1\Big)^{\frac \alp 2}dr+u(x_0)
\\&\geq CR^{2k+\alp}-\widetilde{C},
\end{split}
\]
where $C$ and $\widetilde{C}$ are independent of $R$. Thus there exists a $r_2>2r_1$ large enough such that
\begin{equation}\label{2.13}
\int_R^{+\infty}\frac{R^\alp |\bar{u}(r)|}{r(r^2-R^2)^{\frac\alp 2}}\geq CR^{2k+\alp},\ \ \ \ \forall R>r_2.
\end{equation}
Due to $u\in \mathcal{L}_{\alp}(\R^n)$, we obtain
\begin{equation}\label{2.14}
\int_1^{+\infty}\frac{\overline{|u|}(r)}{r^{1+\alp}}dr=C\int_{|x-x_0|>1}\frac{|u(x)|}{|x-x_0|^{n+\alp}}dx<+\infty.
\end{equation}
Therefore for any $\delta>0$,
\[
\begin{split}
\int_1^{+\infty}&\frac{1}{R^{1+\alp+\delta}}\int_R^{+\infty}\frac{R^\alp |\overline{u}(r)|}{r(r^2-R^2)^{\frac\alp 2}}drdR=\int_1^{+\infty}\frac{|\overline{u}(r)|}{r}\int_1^r\frac{1}{R^{1+\delta}(r^2-R^2)^{\frac\alp 2}}dRdr
\\&\leq C\int_1^{+\infty}\frac{|\overline{u}(r)|}{r^{1+\alp}}\int_1^{\frac r2}\frac{1}{R^{1+\delta}}dRdr+C
\int_1^{+\infty}\frac{|\overline{u}(r)|}{r^{2+\delta}}\int_{\frac r 2}^r\frac{1}{r^{\frac \alp 2}(r-R)^{\frac\alp 2}}dRdr
\\&\leq C_\delta\int_1^{+\infty}\frac{|\overline{u}(r)|}{r^{1+\alp}}dr+C\int_1^{+\infty}\frac{|\overline{u}(r)|}{r^{1+\alp+\delta}}dr
\\&\leq C\int_1^{+\infty}\frac{{|\bar{u}}(r)|}{r^{1+\alp}}dr
\\&\leq C\int_1^{+\infty}\frac{\overline{|u|}(r)}{r^{1+\alp}}dr<+\infty,
\end{split}
\]
which leads to a contradiction with (\ref{2.13}). Thus $k=0$.

From $k=0$, we conclude that
\begin{equation}\label{2.15}
\bar{v}_0(r)\leq\bar{v}_0(0):=-c_0<0, \ \ \ \ \forall r>0.
\end{equation}
Thus for any $R>0$, (\ref{2.4}), (\ref{2.5}), (\ref{2.12}) and (\ref{2.15}) yield that
$$
|\partial B_1|\int_R^{+\infty}\frac{R^\alp| \bar{u}(r)|}{r(r^2-R^2)^{\frac\alp 2}}dr\geq C\int_{0}^{\frac R 2}r^{\alp-1}dr+C\int_{\frac R 2}^R r^{\alp-1}\Big(\frac{R^2}{r^2}-1\Big)^{\frac \alp 2}dr+u(x_0)
\geq CR^{\alp}-\tilde{C},
$$
where positive constants $C$ and $\tilde{C}$ are independent of $R$. Thus there exists $R_0$ large enough such that
\begin{equation}\label{2.16}
\int_R^{+\infty}\frac{R^\alp| \bar{u}(r)|}{r(r^2-R^2)^{\frac\alp 2}}dr\geq CR^\alp,\ \ \ \ \forall R>R_0.
\end{equation}
Obviously
$$\int_N^{+\infty}\frac{\overline{|u|}(r)}{r^{1+\alp}}=C\int_{|x-x_0|>N}\frac{|u(x)|}{|x-x_0|^{n+\alp}}dx=\mbox{o}_N(1)$$
as $N\to+\infty$ due to $u\in \mathcal{L}_{\alp}(\R^n)$. Therefore
\begin{equation}\label{2.17}
\int_{2R}^{+\infty}\frac{R^\alp |\bar{u}(r)|}{r(r^2-R^2)^{\frac\alp 2}}dr\leq CR^{\alp}\int_{2R}^{+\infty}\frac{{|\bar{u}}(r)|}{r^{1+\alp}}\leq \mbox{o}_{R}(1)R^{\alp},
\end{equation}
as $R\to +\infty$. It follows from (\ref{2.16}) and (\ref{2.17}) that there exists $R_1\geq R_0$ large enough such that
\begin{equation}\label{2.18}
\int_R^{2R}\frac{R^\alp| \bar{u}(r)|}{r(r^2-R^2)^{\frac\alp 2}}dr\geq \frac{C}{2}R^\alp,\ \ \ \ \forall R>R_1.
\end{equation}
On the other hand, from (\ref{2.14}), we derive that
\[
\begin{split}
\int_1^{+\infty}\frac{1}{R^{1+\alp}}\int_R^{2R}\frac{R^\alp |\bar{u}(r)|}{r(r^2-R^2)^{\frac\alp 2}}drdR
&\leq \int_1^{+\infty}\frac{|\bar{u}(r)|}{r}\int_{\frac r2}^r\frac{1}{R(r^2-R^2)^{\frac\alp 2}}dRdr
\\&\leq C\int_1^{+\infty}\frac{|\bar{u}(r)|}{r^{2+\frac \alp 2}}\int_{\frac r2}^r\frac{1}{(r-R)^{\frac\alp 2}}dRdr\\&\leq C
\int_1^{+\infty}\frac{{|\bar{u}}(r)|}{r^{1+\alp}}<+\infty,
\end{split}
\]
which contradicts with (\ref{2.18}). Hence Theorem \ref{Th1} is proved.\endproof

\section{Proof of Theorem \ref{T2}}
To complete the proof of Theorem \ref{T2}, we need to establish the following  lemmas.
\begin{lemma}\label{L1}
 Under the  assumptions of Theorem \ref{Th1}, suppose u is a nonconstant classical solution to (1.1) satisfying $u(x)=o(|x|^{\alp+\frac 4 \alp})$ at infinity. Then $(-\Delta)^{\frac \alp 2}u$ satisfies the following integral equation
\begin{align}\label{e1}
(-\Delta)^{\frac \alp 2}u(x)=\int_{\R^n}\frac{R_{2p,n}}{|x-y|^{n-2p}}u_+^\gamma (y)dy+C_p,
\end{align}
where $C_p$ is a nonnegative constant.
\end{lemma}
\proof  Firstly, we prove that
\begin{equation}\label{3.1}
(-\Delta)^{p-1+\frac \alp 2}u(x)=\int_{\R^n}\frac{R_{2,n}}{|x-y|^{n-2}}u_+^\gamma (y)dy,\ \ \ \ \forall x\in \R^n.
\end{equation}

To this end, for arbitrary $R>0$, denote $f_1(x)=u_+^\gamma(x)$ and
$$v_1^R(x):=\int_{B_R(x_0)}G_R^2(x,y)f_1(y)dy,$$
where the Green's function for $-\Delta$ on $B_R(0)$ is given by
$$G_R^2(x,y)=R_{2,n}\Big[\frac{1}{|x-y|^{n-2}}-\frac{1}{(|x|\cdot|\frac{Rx}{|x|^2}-\frac{y}{R}|)^{n-2}}\Big], \ \ \ \ x, y \in B_R(0),$$
and $G_R^2(x,y)=0$ if $x$ or $y\in\R^n\setminus B_R(0)$. Then we derive that $(-\Delta)^{p-1+\frac \alp 2}u-v_1^R$ is harmonic in $B_R(0)$ and continuous up to $\partial B_R(0)$. Hence $v_1^R(x)\in C^2(B_R(0))\cap C(\R^n)$ and satisfies
\begin{equation}\label{3.2}
\left\{\begin{aligned}
&-\Delta v_1^R(x)=u_+^\gamma,&\ \ \ \ &x\in B_R(0) ,
\\& v_1^R(x)=0,&\ \ \ \ &x\in\R^n\setminus B_R(0).
\end{aligned}\right.
\end{equation}
Set $w_1^R(x):=(-\Delta)^{p-1+\frac \alp 2}u-v_1^R$. Combining Theorem \ref{Th1} with (\ref{3.2}), we obtain $w_1^R(x)\in C^2(B_R(0))\cap C(\R^n)$ and satisfies
$$
\left\{\begin{aligned}
&-\Delta w_1^R(x)=0,&\ \ \ \ &x\in B_R(0) ,
\\& w_1^R(x)\geq0,&\ \ \ \ &x\in\R^n\setminus B_R(0).
\end{aligned}\right.
$$
The maximum principle implies that
\begin{equation}\label{3.3}
w_1^R(x)=(-\Delta)^{p-1+\frac \alp 2}u(x)-v_1^R(x)\geq 0,\ \ \ \ \forall x \in \R^n.
\end{equation}
Now for each fixed $x\in \R^n$, letting $R\to +\infty$ in (\ref{3.3}), we have
\begin{equation}\label{3.4}
(-\Delta)^{p-1+\frac \alp 2}u(x)\geq\int_{\R^n}\frac{R_{2,n}}{|x-y|^{n-2}}f_1(y)dy=:v_1(x)\geq 0.
\end{equation}
Since $v_1-v_1^R$ is smooth in $B_R(0)$ for any $R>0$, one can obtain that $v_1\in C^2(\R^n)$ satisfies
\begin{equation}\label{3.5}
-\Delta v_1(x)=u_+^ \gamma(x), \ \ \ \  \forall x\in \R^n.
\end{equation}
Set $w_1(x):=(-\Delta)^{p-1+\frac \alp 2}u(x)-v_1(x)$. From (\ref{T1}), (\ref{3.4}) and (\ref{3.5}), we have that $w_1(x)\in C^2(\R^n)$ and satisfies
$$
\left\{\begin{aligned}
&-\Delta w_1(x)=0,&\ \ \ \ &x\in \R^n ,
\\& w_1(x)\geq0,&\ \ \ \ &x\in\R^n.
\end{aligned}\right.
$$
Then by the Liouville Theorem for harmonic function, we can infer that
$$w_1(x)=(-\Delta)^{p-1+\frac \alp 2}u(x)-v_1(x)\equiv C_1\geq 0.$$
Therefore
\begin{equation}\label{3.6}
(-\Delta)^{p-1+\frac \alp 2}u(x)=\int_{\R^n}\frac{R_{2,n}}{|x-y|^{n-2}}u_+^\gamma (y)dy+C_1=:f_2(x)\geq C_1\geq 0.
\end{equation}

Next, for arbitrary $R>0$, denote
$$v_2^R(x):=\int_{B_R(x_0)}G_R^2(x,y)f_2(y)dy.$$
Then we have
\begin{equation}\label{3.7}
\left\{\begin{aligned}
&-\Delta v_2^R(x)=f_2(x),&\ \ \ \ &x\in B_R(0) ,
\\& v_2^R(x)=0,&\ \ \ \ &x\in\R^n\setminus B_R(0).
\end{aligned}\right.
\end{equation}
Set $w_2^R(x):=(-\Delta)^{p-2+\frac \alp 2}u-v_2^R$. By Theorem \ref{Th1}, (\ref{3.6}) and (\ref{3.7}), we obtain
$$
\left\{\begin{aligned}
&-\Delta w_2^R(x)=0,&\ \ \ \ &x\in B_R(0) ,
\\& w_2^R(x)\geq0,&\ \ \ \ &x\in\R^n\setminus B_R(0).
\end{aligned}\right.
$$
The maximum principle implies that for any $R>0$
\begin{equation}\label{3.8}
w_2^R(x)=(-\Delta)^{p-2+\frac \alp 2}u(x)-v_2^R(x)\geq 0,\ \ \ \ \forall x \in \R^n.
\end{equation}
Now for each fixed $x\in \R^n$, letting $R\to +\infty$ in (\ref{3.8}), we have
\begin{equation}\label{3.9}
(-\Delta)^{p-2+\frac \alp 2}u(x)\geq\int_{\R^n}\frac{R_{2,n}}{|x-y|^{n-2}}f_2(y)dy=:v_2(x)\geq 0.
\end{equation}
Taking $x=0$ in (\ref{3.9}) and using (\ref{3.6}), we get
$$\int_{\R^n}\frac{C_1}{|y|^{n-2}}dy\leq\int_{\R^n}\frac{f_2(y)}{|y|^{n-2}}dy<+\infty.$$
It follows immediately that $C_1=0$, and hence we shows that (\ref{3.1}) holds, that is
\begin{equation}\label{3.10}
(-\Delta)^{p-1+\frac \alp 2}u(x)=f_2(y)=\int_{\R^n}\frac{R_{2,n}}{|x-y|^{n-2}}u_+^\gamma (y)dy,\ \ \ \ \forall x\in \R^n.
\end{equation}

One can easily observe that $v_2$ is a solution of
\begin{equation}\label{3.11}
-\Delta v_2(x)=f_2(x), \ \ \ \  \forall x\in \R^n.
\end{equation}
Set $w_2(x):=(-\Delta)^{p-2+\frac \alp 2}u(x)-v_2(x)$. Then it solves
$$
\left\{\begin{aligned}
&-\Delta w_2(x)=0,&\ \ \ \ &x\in \R^n ,
\\& w_2(x)\geq0,&\ \ \ \ &x\in\R^n.
\end{aligned}\right.
$$
Applying the Liouville Theorem for harmonic functions again, we can infer that
$$w_2(x)=(-\Delta)^{p-2+\frac \alp 2}u(x)-v_2(x)\equiv C_2\geq 0.$$
Therefore, we deduce that
\begin{equation}\label{3.12}
(-\Delta)^{p-2+\frac \alp 2}u(x)=\int_{\R^n}\frac{R_{2,n}}{|x-y|^{n-2}}f_2 (y)dy+C_2=:f_3(x)\geq C_2\geq 0.
\end{equation}
By the same methods as above, we can prove that $C_2=0$, and hence
$$(-\Delta)^{p-2+\frac \alp 2}u(x)=f_3(y)=\int_{\R^n}\frac{R_{2,n}}{|x-y|^{n-2}}f_2 (y)dy,\ \ \ \ \forall x\in \R^n.$$
Define
\begin{equation}\label{3.13}
f_{k+1}(x):=\int_{\R^n}\frac{R_{2,n}}{|x-y|^{n-2}}f_k (y)dy
\end{equation}
for $k=1,2,\ldots,p$. Repeating the above argument we can derive that
\begin{equation}\label{3.14}
(-\Delta)^{p-k+\frac \alp 2}u(x)=\int_{\R^n}\frac{R_{2,n}}{|x-y|^{n-2}}f_k (y)dy,\ \ \ \ \forall x\in \R^n
\end{equation}
for $k=1,2,\ldots,p-1$ and
\begin{equation}\label{3.15}
(-\Delta)^{\frac \alp 2}u(x)=\int_{\R^n}\frac{R_{2,n}}{|x-y|^{n-2}}f_p (y)dy+C_p=f_{p+1}(x)+C_p\geq C_p\geq 0,\ \ \ \ \forall x\in \R^n.
\end{equation}
Owing to the properties of the Riesz potential, for any $\alp_1, \alp_2\in(0,n)$ such that $\alp_1+\alp_2\in(0,n)$, it is easy to know that
\begin{equation}\label{3.16}
\int_{\R^n}\frac{R_{\alp_1,n}}{|x-y|^{n-\alp_1}}\cdot\frac{R_{\alp_2,n}}{|y-z|^{n-\alp_2}}dy
=\frac{R_{\alp_1+\alp_2,n}}{|x-z|^{n-(\alp_1+\alp_2)}}.
\end{equation}
In particular, it follows from (\ref{3.13})-(\ref{3.16}) and Fubini's theorem that
$$(-\Delta)^{p-k+\frac \alp 2}u(x)=\int_{\R^{n}}\frac{R_{2k,n}}{|x-y|^{n-2k}}u_+^\gamma (y)dy,\ \ \ \ \forall x\in \R^n$$
for $k=1,2,\ldots,p-1$ and
\begin{equation*}
(-\Delta)^{\frac \alp 2}u(x)=\int_{\R^n}\frac{R_{2p,n}}{|x-y|^{n-2p}}u_+^\gamma (y)dy+C_p,\ \ \ \ \forall x\in \R^n.
\end{equation*}
Thus the lemma holds.\endproof

\begin{lemma}\label{L2}
Under the  assumptions of Theorem \ref{Th1}, if $ u\in C_{loc}^{2p+[\alp],\{\alp\}+\eps}(\R^n)$ satisfies $u_+^\gamma\in L^1(\R^n)$ and (\ref{e1}), then $\Delta u\in L^\infty(\R^n)$.
\end{lemma}
\proof From (\ref{e1}), (\ref{3.16}) and the formula of fundamental solution for $(-\Delta)^{\frac \alp 2}$, we obtain
\begin{align}\label{a1}
-\Delta u(x)&=(-\Delta)^{1-\frac \alp 2}\Big(\int_{\R^n}\frac{R_{2p,n}}{|x-y|^{n-2p}}u_+^\gamma (y)dy+C_p\Big)
\notag\\&=(-\Delta)^{1-\frac \alp 2}\Big(\int_{\R^n}u_+^\gamma (y)\int_{\R^n}\frac{R_{2-\alp,n}}{|x-z|^{n-2+\alp}}\frac{R_{2p-2+\alp,n}}{|z-y|^{n-2p+2-\alp}}dzdy\Big)
\notag\\&=(-\Delta)^{1-\frac \alp 2}\Big(\int_{\R^n}\frac{R_{2-\alp,n}}{|x-z|^{n-2+\alp}} \int_{\R^n}\frac{R_{2p-2+\alp,n}}{|z-y|^{n-2p+2-\alp}}u_+^\gamma (y)dydz\Big)
\notag\\&={R_{2p-2+\alp,n}}\int_{\R^n}\frac{u_+^\gamma (y)}{|x-y|^{n-2p+2-\alp}}dy.
\end{align}

Then analogous to Lemma 2.3 in \cite{Du}, we can immediately infer that  $\Delta u\in L^\infty(\R^n)$. Indeed,
from (\ref{a1}), we derive that
\begin{align}\label{limit8}
~\Delta u(x)+R_{2p-2+\alp,n}\int_{\mathbb{R}^n}\frac{u_+^\gamma(y)}{|x-y|^{n-2p+2-\alp}}dy=0.
\end{align}
Denote $\beta:=\int_{\mathbb{R}^n}u_+^\gamma(y)dy$.  For any $x_0\in{\mathbb{R}^n}$, consider the solution $h_1$ of the boundary value problem
\begin{align}\label{limit17}
\left\{\begin{array}{ll}
-\Delta h_1(x)=R_{2p-2+\alp,n}\int_{B_{{R}}(x_0)}\frac{u_+^\gamma(y)}{|x-y|^{n-2p+2-\alp}}dy  ~~ &\mbox{ in }{B_{{R}}(x_0)},\\
h_1=0~~~~~~~~~~~~~~~~~~~~~~~~~~~~~~~~&\mbox{ on } {\partial B_{{R}}(x_0)},
\end{array}\right.
\end{align}
where ${R}$ will be determined later. Set
$$ v_1(x)=C_{n,p,\alp}\int_{B_{{R}}(x_0)}\frac{u_+^\gamma(y)}{|x-y|^{n-2p-\alp}}dy ~~\mbox{ for any } x\in B_{{R}}(x_0),$$
where $C_{n,p,\alp}=\frac{R_{2p-2+\alp,n}}{(n-2p-\alp)(2p-2+\alp)}$. It's obvious that $v_1(x)\geq 0$ in $B_{{R}}(x_0)$ and
\begin{align}\label{limit10}
-\Delta v_1(x)=R_{2p-2+\alp,n}\int_{B_{{R}}(x_0)}\frac{u_+^\gamma(y)}{|x-y|^{n-2p+2-\alp}}dy.
\end{align}
The comparison  principle allows us to conclude that
$$ |h_1(x)|\leq v_1(x), ~~x\in B_{{R}}(x_0).$$
Therefore there exists $t>1$ and a constant $C>0$ independent of $x_0$ such that
\begin{align}\label{limit11}
\int_{B_{{R}}(x_0)}|h_1(x)|^{t\gamma}dx\leq C.
\end{align}
Indeed
\begin{align}\label{limit15}
\begin{split}
\int_{B_{{R}}(x_0)}|h_1(x)|^{t\gamma}dx
\leq\int_{B_{{R}}(x_0)}|v_1(x)|^{t\gamma}dx
=\int_{B_{{R}}(x_0)}\left|{C}_{n,p,\alp}\int_{B_{{R}}(x_0)}\frac{u_+^\gamma(y)}{|x-y|^{n-2p-\alp}}dy\right|^{t\gamma}dx.
\end{split}
\end{align}
Owing to $\gamma\in(1,\frac{n}{n-2p-\alp})$, we know that there exists $t>1$ such that $t\gamma <\frac{n}{n-2p-\alp}$. Denote $d\mu=u_+^\gamma dy/ \int_{B_{{R}}(x_0)}u_+^\gamma(y)dy$, then the assumption $u_+^\gamma(y)\in L^1(\mathbb{R}^n)$ and Jensen's inequality imply  that
\begin{align}\label{limit16}
\begin{split}
\int_{B_{{R}}(x_0)}\left|{\int_{B_{{R}}(x_0)}\frac{u_+^\gamma(y)}{|x-y|^{n-2p-\alp}}dy}\right|^{t\gamma}dx
\leq& C\int_{B_{{R}}(x_0)}\left(\int_{B_{{R}}(x_0)}\frac{1}{|x-y|^{n-2p-\alp}}d\mu\right)^{t\gamma}dx\\
\leq&  C\int_{B_{{R}}(x_0)}\int_{B_{{R}}(x_0)}\frac{1}{|x-y|^{(n-2p-\alp)t\gamma}}dxd\mu \\
\leq& C,~~~~
\end{split}
\end{align}
where the final inequality used the fact that $t\gamma(n-2p-\alp)<n$.
From this and (\ref{limit15}), we  obtain  (\ref{limit11}).

Now, we consider the function $q(x):=u(x)-h_1(x)$ in the smaller ball $B_{{R}-1}(x_0)$. First we observe that
\begin{align*}
\begin{split}
\Delta q(x)&=\Delta u(x)-\Delta h_1(x)=-R_{2p-2+\alp,n}\int_{\mathbb{R}^n\setminus B_{{R}}(x_0)}\frac{u_+^\gamma(y)}{|x-y|^{n-2p+2-\alp}}dy.
\end{split}
\end{align*}
If $x\in B_{{R}-1}(x_0)$ and $y\in \mathbb{R}^n\setminus B_{{R}}(x_0)$, then $|x-y|\geq 1$. Therefore we have
\begin{align}\label{limit12}
0\leq -\Delta q(x)\leq R_{2p-2+\alp,n} \beta.
\end{align}
Hence, it follows from weak Harnack principle that
\begin{align}\label{limit13}
\sup\limits_{B_{{R}-2}(x_0)}q(x)\leq C[\|q_+\|_{L^\gamma(B_{{R}-1}(x_0))}+\|\Delta q\|_{L^\infty(B_{{R}-1}(x_0))}].
\end{align}
 (\ref{limit12}) shows that the second term of the right hand side of the above inequality is bounded independent of $x_0$. By using (\ref{limit11}) and $u_+^\gamma\in L^1(\mathbb{R}^n)$, we obtain that
$$ \int_{B_{{R}-1}(x_0)}(q_+(x))^\gamma dx\leq C,$$
where $C$ independent of $x_0$. Therefore, it follows that $u(y)=q(y)+h_1(y)\leq C+|h_1(y)|$ in the smaller ball $B_{{R}-2}(x_0)$. This and (\ref{limit11}) yield the estimate
\begin{align}\label{limit19}
\int_{B_{{R}-2}(x_0)}u_+^{t\gamma}(y)dy\leq C.
\end{align}

Now consider the solution $h_2$ of the boundary value problem
\begin{align}\label{limit2338}
\left\{\begin{array}{ll}
-\Delta h_2(x)=R_{2p-2+\alp,n}\int_{B_{{R}-2}(x_0)}\frac{u_+^\gamma(y)}{|x-y|^{n-2p+2-\alp}}dy  ~~ &\mbox{ in }{B_{{R}-2}(x_0)},\\
h_2=0~~~~~~~~~~~~~~~~~~~~~~~~~~~~~~~~~~&\mbox{ on } {\partial B_{{R}-2}(x_0)}.
\end{array}\right.
\end{align}
Set
$$ v_2(x)=C_{n,p,\alp}\int_{B_{{R}-2}(x_0)}\frac{u_+^\gamma(y)}{|x-y|^{n-2p-\alp}}dy ~~\mbox{ for any } x\in B_{{R}-2}(x_0),$$
and we have that $|h_2(x)|\leq v_2(x)$ in $B_{{R}-2}(x_0)$. Simple computation shows that
\begin{eqnarray}\nonumber
\int_{B_{{R}-2}(x_0)}\frac{u_+^\gamma(y)}{|x-y|^{n-2p-\alp}}dy&&\leq C \left(\int_{B_{{R}-2}(x_0)}\frac{u_+^{t\gamma}(y)}{|x-y|^{n-2p-\alp}}dy\right)^{\frac{1}{t}}
\left(\int_{B_{{R}-2}(x_0)}\frac{1}{|x-y|^{n-2p-\alp}}dy\right)^{\frac{t-1}{t}}
\\ \nonumber
&&\leq C \left(\int_{B_{{R}-2}(x_0)}\frac{u_+^{t\gamma}(y)}{|x-y|^{n-2p-\alp}}dy\right)^{\frac{1}{t}}.
\end{eqnarray}
Hence
\begin{align}\label{limit25}
\begin{split}
\int_{B_{{R}-2}(x_0)}|h_2(x)|^{t^2\gamma}dx &\leq \int_{B_{{R}-2}(x_0)}|v_2(x)|^{t^2\gamma}dx \\
&= C \int_{B_{{R}-2}(x_0)}\left(\int_{B_{{R}-2}(x_0)}\frac{u_+^\gamma(y)}{|x-y|^{n-2p-\alp}}dy\right)^{t^2\gamma}dx\\
&\leq C \int_{B_{{R}-2}(x_0)}\left(\int_{B_{{R}-2}(x_0)}\frac{u_+^{t\gamma}(y)}{|x-y|^{n-2p-\alp}}dy\right)^{t\gamma}dx.
\end{split}
\end{align}
Let $d\mu=u_+^{t\gamma} dy/ \int_{B_{{R}-2}(x_0)}u_+^{t\gamma}(y)dy$.
Therefore (\ref{limit19}), (\ref{limit25}) and Jensen's inequality imply  that there exists a constant $C>0$ independent of $x_0$ such that
\begin{align*}
\begin{split}
\int_{B_{{R}-2}(x_0)}|h_2(x)|^{t^2\gamma}dx \leq C.
\end{split}
\end{align*}
Similar argument as the above, we can obtain that there exists a constant $C>0$ independent of $x_0$ such that
\begin{align}\label{limit07}
\int_{B_{{R}-4}(x_0)}u_+^{t^2\gamma}(y)dy\leq C.
\end{align}

Let $h_3$ be the solution of the equation
\begin{align}\label{limit2352}
\left\{\begin{array}{ll}
-\Delta h_3(x)=R_{2p-2+\alp,n}\int_{B_{{R}-4}(x_0)}\frac{u_+^\gamma(y)}{|x-y|^{n-2p+2-\alp}}dy  ~~ &\mbox{ in }{B_{{R}-4}(x_0)},\\
h_3=0,~~~~~~~~~~~~~~~~~~~~~~~~~~~~~~~~~~&\mbox{ on } {\partial B_{{R}-4}(x_0)},
\end{array}\right.
\end{align}
and
$$ v_3(x)=C_{n,p,\alp}\int_{B_{{R}-4}(x_0)}\frac{u_+^\gamma(y)}{|x-y|^{n-2p-\alp}}dy ~~\mbox{ for all } x\in B_{{R}-4}(x_0).$$
Similar argument as before yields that
$$ \int_{B_{{R}-6}(x_0)}u_+^{t^3\gamma}(y)dy\leq C.$$
Repeating the process, we  derive that for any $k\in \mathbb{N}$ satisfying $R-2k>0$,
\begin{align}\label{limit22}
\int_{B_{{R}-2k}(x_0)}u_+^{t^k\gamma}(y)dy\leq C.
\end{align}

Choose $k>0$ large enough such that $\frac{1}{t^k}+\frac{1}{s}=1$ where $s$ satisfies $s(n-2p+2-\alp)<n$. For such $k$, we choose ${R}$ large enough such that ${{R}-2k}>2$. From (\ref{limit8}), we have
\begin{align*}
\begin{split}
|\Delta u(x_0)|
&=R_{2p-2+\alp,n}\int_{\mathbb{R}^n\setminus {B_{{R}-2k}(x_0)}}\frac{u_+^\gamma(y)}{|x_0-y|^{n-2p+2-\alp}}dy+R_{2p-2+\alp,n}\int_{{B_{{R}-2k}(x_0)}}\frac{u_+^\gamma(y)}{|x_0-y|^{n-2p+2-\alp}}dy\\
&\leq \frac{R_{2p-2+\alp,n}}{4}\int_{\mathbb{R}^n}u_+^\gamma (y)dy +C\left(\int_{{B_{{R}-2k}(x_0)}}u_+^{t^k\gamma}(y)dy\right)^{\frac{1}{t^k}}
\left(\int_{{B_{{R}-2k}(x_0)}}\frac{1}{|x_0-y|^{(n-2p+2-\alp)s}}dy\right)^{\frac{1}{s}}\\
&\leq \frac{R_{2p-2+\alp,n}}{4}\int_{\mathbb{R}^n}u_+^\gamma (y)dy +C\left(\int_{{B_{{R}-2k}(x_0)}}u_+^{t^k\gamma}(y)dy\right)^{\frac{1}{t^k}},
\end{split}
\end{align*}
where the final inequality used the fact that $(n-2p+2-\alp)s<n$.
(\ref{limit22}) and the integral constraint in (\ref{T1}) show that $|\Delta u(x_0)|\leq C$, where $C>0$ is independent of $x_0$. Therefore, we finish the proof of this lemma.
$\hfill\square$

\begin{lemma}\label{a2}
If $\gamma\in(1,\frac{n}{n-2p-\alp})$, $n>2p+\alp$, $p\geq1$ is an integer and $u$ is a solution of (\ref{T1}) satisfying $u(x)=o(|x|^{\alp +\frac 4 \alp})$ at infinity. then there exists a constant $M>0$ such that $\sup\limits_{\mathbb{R}^n}u \leq M$.
\end{lemma}
{\bf Proof} From Lemma  \ref{L2}, we have that there exists $A>0$ such that $|\Delta u| \leq A$ in ${\mathbb{R}^n}.$
Denote $h(x)=-\Delta u(x)$. Given
$x_0\in {\mathbb{R}^n}$, let $u_1$ be the solution of
\begin{align}\label{03}
\left\{\begin{array}{rl}
-\Delta v=h, ~~~&\mbox{ in } B_1(x_0),\\
v=0,~~~~~~~~~~ &\mbox{ on } \partial B_1(x_0).
\end{array}\right.
\end{align}
It follows from the elliptic theory that $|u_1|\leq C$, where $C>0$ independent of $x_0$. Denote $u_2=u-u_1$,
then $(u_2)_+\leq u_++|u_1|$. Since $|u_1|\leq C$ in $B_1(x_0)$ and $\int_{\mathbb{R}^n}u_+^{\gamma}(x)dx<+\infty$, we derive
$$\int_{B_1(x_0)}(u_2)_+^{\gamma}(x)\leq C.$$
Note that $\Delta u_2=0$ in $B_1(x_0)$.
For the subharmonic function $(u_2)_+$, we have
$$\|(u_2)_+\|_{L^\infty(B_{1/2}(x_0))}\leq C \int _{B_1(x_0)}(u_2)_+(x)dx \leq C \left(\int _{B_1(x_0)}(u_2)_+^\gamma(x)dx\right)^{\frac{1}{\gamma}}\leq C.$$
where $C$ is independent of $x_0$.
Recalling that $u=u_1+u_2$ and the arbitrariness of $x_0$, we derive that there exists $M>0$ independent of $x_0$ such that $u_+(x)\leq M$. Hence, $\sup\limits_{\mathbb{R}^n}u \leq M$.
$\hfill\square$

Set
$$\zeta(x)=-\int_{\R^n}\frac{R_{2p+\alp,n}}{|x-y|^{n-2p-\alp}}u_+^\gamma (y)dy.$$
Then it is easy to obtain some asymptotic behaviors about $\zeta(x)$ at infinity as follows.
\begin{lemma}\label{L4}
$\zeta(x)$ satisfies
$$\lim_{|x|\to+\infty}(-\Delta)^i\zeta(x)|x|^{n-2p-\alp+2i}=a_i, \ \ \ \ i=0,1,\ldots,p-1,$$
where $a_0=-R_{2p+\alp,n}\int_{\R^n}u_+^\gamma(y)dy$ and $a_{i+1}=a_i(n-2p-\alp+2i)(2p+\alp-2i-2)$, $i=0,1,\ldots, p-2$.

Moreover,
$$\lim_{|x|\to+\infty}(-\Delta)^{i+\frac{\alp}{2}}\zeta(x)|x|^{n-2p+2i}=b_i, \ \ \ \ i=0,1,\ldots,p-1,$$
where $b_0=-R_{2p,n}\int_{\R^n}u_+^\gamma(y)dy$ and $b_{i+1}=b_i(n-2p+2i)(2p-2i-2)$, $i=0,1,\ldots, p-2$.
\end{lemma}

We are now ready to complete the proof of the equivalence between (\ref{T1}) and (\ref{1.2}), that is proof of Theorem \ref{T2}.

{\bf Proof of Theorem \ref{T2}}
 From Lemma {\ref{L1}}, we derive (\ref{e1}) holds. Next we will show that $C_p=0$.

For arbitrary $R>0$, let
$$v_{p+1}^R(x):=\int_{B_R(x_0)}G_R^\alp(x,y)(f_{p+1}(y)+C_p)dy,$$
where $G_R^\alp(x,y)$ is the Green's function for $(-\Delta)^{\frac \alp 2}$ with $0<\alp<2$ on $B_R(0).$
Then, we can get
\begin{equation}\label{3.18}
\left\{\begin{aligned}
&(-\Delta)^{\frac \alp 2} v_{p+1}^R(x)=f_{p+1}(x)+C_p,&\ \ \ \ &x\in B_R(0) ,
\\& v_{p+1}^R(x)=0,&\ \ \ \ &x\in\R^n\setminus B_R(0).
\end{aligned}\right.
\end{equation}
Denote $w_{p+1}^R(x):=M-u(x)+v_{p+1}^R(x)$. From (\ref{e1}), (\ref{3.18}) and Lemma \ref{a2}, we have
$$
\left\{\begin{aligned}
&(-\Delta)^{\frac \alp 2} w_{p+1}^R(x)=0,&\ \ \ \ &x\in B_R(0) ,
\\& w_{p+1}^R(x)\geq0,&\ \ \ \ &x\in\R^n\setminus B_R(0).
\end{aligned}\right.
$$
By maximum principle, we can deduce that for any $R>0$
\begin{equation}\label{3.19}
w_{p+1}^R(x)=M-u(x)+v_{p+1}^R(x)\geq 0,\ \ \ \ \forall x \in \R^n.
\end{equation}
Now for each fixed $x\in \R^n$, letting $R\to +\infty$ in (\ref{3.19}), we have
\begin{equation}
M-u(x)\geq-\int_{\R^n}\frac{R_{\alp,n}}{|x-y|^{n-\alp}}(f_{p+1}(y)+C_p)dy=:-v_{p+1}(x).
\end{equation}
Thus as $R\to +\infty$, we have
$$
\left\{\begin{aligned}
&(-\Delta)^{\frac \alp 2}(M-u(x)+v_{p+1}(x))=0,&\ \ \ \ &x\in B_R(0) ,
\\& M-u(x)+v_{p+1}(x)\geq0,&\ \ \ \ &x\in\R^n\setminus B_R(0).
\end{aligned}\right.
$$
By Liouville Theorem, we obtain
$$ M-u(x)+v_{p+1}(x)\equiv \widetilde{C}\geq0.$$ That is,
\begin{equation}\label{3.21}
u(x)=\int_{\R^n}\frac{R_{\alp,n}}{|x-y|^{n-\alp}}(f_{p+1}(y)+C_p)dy+C.
\end{equation}
Taking $x=0$ in (\ref{3.21}), we get
$$\int_{\R^n}\frac{C_p}{|y|^{n-\alp}}dy\leq\int_{\R^n}\frac{f_{p+1}(y)+C_p}{|y|^{n-\alp}}dy<+\infty.$$
It follows immediately that $C_p=0$. Thus, from (\ref{3.13}), (\ref{3.16}) and (\ref{3.21}), we derive that (\ref{1.2}) holds. Moreover, we assert that $C_0<0$. Indeed, if  $C_0\geq0$, we have $u(x)\geq 0$ in $\R^n$, which is impossible from $\gamma \in (1,\frac{n+2p+\alp}{n-2p-\alp})$ and the results of Theorem 1.9 in \cite{3}. Hence (\ref{1.2}) holds.

From Lemma \ref{L4} and the fact $C_0<0$, we  obtain  that the support of $u_+$ is compact. Meanwhile, it is obvious that if $u$ is a solution of (\ref{1.2}), then it satisfies equation (\ref{T1}). The proof of Theorem \ref{T2} is  completed.\endproof

\section{Proof of Theorem 3}
In this section, we verify Theorem \ref{T3} by taking advantage of the method of moving planes in integral forms.

{\bf Proof }  To complete the proof of Theorem \ref{T3}, it's enough to show that $\zeta$ is symmetric about some point $x_0\in{\mathbb{R}^n}$ and $\frac{\partial \zeta}{\partial r}>0$ where $r=|x-x_0|$. From Lemmas \ref{L2}-\ref{L4}, we obtain that $\lim\limits_{|x|\rightarrow \infty}\zeta(x)=0$ and
\begin{equation}\label{x2}
\zeta(x)=-R_{2p+\alp,n}\int_{\mathbb{R}^n}\frac{(C_0-\zeta)_+^\gamma(y)}{|x-y|^{n-2p-\alp}}dy, ~~x\in{\mathbb{R}^n}.
\end{equation}
For $x=(x_1,x_2,...,x_n)\in {\mathbb{R}}^n$ and $\lambda\in{\mathbb{R}}$, we define ${\mathrm{T}}_{\lambda}=\{x\in{\mathbb{R}}^n|x_1=\lambda\}$, $\Sigma_\lambda=\{x\in{\mathbb{R}}^n|x_1<\lambda\}$, $x^\lambda=(2\lambda-x_1,x_2,...,x_n)$ and $\zeta_\lambda(x)=\zeta(2\lambda-x_1,x_2,...,x_n)=\zeta(x^\lambda)$.
Set $w_\lambda(x)=\zeta(x)-\zeta_\lambda(x)$. It's obvious that
\begin{equation*}
\zeta_\lambda(x)=-R_{2p+\alp,n}\int_{\mathbb{R}^n}\frac{(C_0-\zeta_\lambda)_+^\gamma(y)}{|x-y|^{n-2p-\alp}}dy, ~~x\in{\mathbb{R}^n}.
\end{equation*}
From this and (\ref{x2}), we have
\begin{align}\label{x5}
\zeta_\lambda(x)-\zeta(x)=R_{2p+\alp,n}\int_{\Sigma_\lambda}\left(\frac{1}{|x-y|^{n-2p-\alp}}-\frac{1}{|x-y^\lambda|^{n-2p-\alp}}\right)\left((C_0-\zeta)_+^\gamma(y)-(C_0-\zeta_\lambda)_+^\gamma(y)\right)dy.
\end{align}
\par $\mathbf{Step~ 1}$: We claim that for $\lambda$ sufficiently negative,
\begin{align}\label{x3}
w_\lambda(x)> 0, ~~~~~x\in \Sigma_\lambda.
\end{align}
Due to $\lim\limits_{|x|\rightarrow \infty}\zeta(x)=0$ and $C_0<0$, we have for $\lambda$ sufficiently negative
$$ (C_0-\zeta)_+^\gamma-(C_0 -\zeta_\lambda)_+^\gamma=-(C_0 -\zeta_\lambda)_+^\gamma\leq 0,~~~~x\in \Sigma_\lambda.$$
From this, (\ref{x5}) and the fact that $u$ is a sign-changing classical solution of (\ref{T1}), we have $\zeta_\lambda(x)-\zeta(x)< 0$ for any $x\in \Sigma_\lambda$. Thus (\ref{x3}) holds.

\par $\mathbf{Step~ 2}$: Step 1 provides a starting point, from which we can now move the plane $T_\lambda$ to the right as long as (\ref{x3}) holds to its limiting position. Define
\begin{equation*}
\lambda_0=\sup\{\lambda |w_\mu(x)>0, \forall x \in \Sigma_\mu, \mu\leq\lambda\}.
\end{equation*}
It's obviously that $\lambda_0<+\infty$ and
$$ w_{\lambda_0}(x)\geq 0, ~~~x\in\Sigma_{\lambda_0}.$$
We will show that $ w_{\lambda_0}(x)\equiv0$ for $x\in\Sigma_{\lambda_0}$.

Otherwise if $w_{\lambda_0}\geq 0$ and $w_{\lambda_0}\not\equiv0$, we must have
\begin{align}\label{x7}
w_{\lambda_0}(x)>0,~~~x\in\Sigma_{\lambda_0},
\end{align}
where (\ref{x7}) follows from (\ref{x5}). We can derive that there exists $R$ large enough such that
\begin{equation}\label{p1}
(C_0-\zeta)_+\equiv0, ~~x\in\mathbb{R}^n\setminus B_R(0),
\end{equation}
due to  $\lim\limits_{|x|\rightarrow \infty}\zeta(x)=0$ and $C_0<0$. Fixing this $R$, we have there exists constant $\delta>0$ and $c>0$ such that
\begin{align}\label{x8}
w_{\lambda_0}(x)\geq c,   ~~~x\in\overline{\Sigma_{\lambda_0 -\delta}\cap B_R(0)}.
\end{align}
Therefore by the continuity of $w_\lambda$ in $\lambda$ there exists $\varepsilon>0$ and $\varepsilon<\delta$ such that for all $\lambda\in[\lambda_0,\lambda_0+\varepsilon)$, we have
\begin{equation*}
w_\lambda(x)\geq0, ~~~x\in\overline{\Sigma_{\lambda_0 -\delta}\cap B_R(0)}.
\end{equation*}
We will show that for sufficiently small $0<\varepsilon<\delta$ and any $\lambda\in[\lambda_0,\lambda_0+\varepsilon)$
\begin{align}\label{x11}
w_\lambda(x)\geq0, ~~~x\in\Sigma_\lambda,
\end{align}
which contradicts with the definition of $\lambda_0$. Therefore we must have $w_{\lambda_0}\equiv 0$. Define
$$\Sigma_\lambda^-=\{x\in\Sigma_\lambda|w_\lambda(x)<0\}.$$
Next we claim that $\Sigma_\lambda^-$ must be measure zero.

For $y\in\Sigma_\lambda^-$, we can obtain that
\begin{equation}\label{4.1}
(C_0 -\zeta)_+^\gamma(y)-(C_0-\zeta_\lambda)_+^\gamma(y)\leq \gamma (C_0 -\zeta)_+^{\gamma-1}(y)|w_\lambda(y)|.
\end{equation}
Thus for $x\in \Sigma_\lambda$,
\begin{align}\label{x22}
\begin{split}
\zeta_\lambda(x)-\zeta(x)&\leq R_{2p+\alp,n}\int_{\Sigma_\lambda^-}\left(\frac{1}{|x-y|^{n-2p-\alp}}-\frac{1}{|x-y^\lambda|^{n-2p-\alp}}\right)\left((C_0- \zeta)_+^\gamma(y)-(C_0-\zeta_\lambda)_+^\gamma(y)\right)dy \\
&\leq R_{2p+\alp,n}\int_{\Sigma_\lambda^-}\left(\frac{1}{|x-y|^{n-2p-\alp}}-\frac{1}{|x-y^\lambda|^{n-2p-\alp}}\right)\gamma (C_0-\zeta(y))_+^{\gamma-1}|w_\lambda(y)|dy.
\end{split}
\end{align}
Applying Hardy-Littlewood-Sobolev inequality \cite{Dai} and H\"{o}lder inequality to (\ref{x22}) we obtain that
\begin{align}\label{x21}
\begin{split}
\|\zeta_\lambda(x)-\zeta(x)\|_{L^{\frac{2n}{n-2p-\alp}}(\Sigma_\lambda^-)}&\leq C\left(\int_{\Sigma_\lambda^-}\left((C_0-\zeta)_+^{\gamma-1}(y)|w_{\lambda}(y)|\right)^{\frac{2n}{n+2p+\alp}}dy\right)^{\frac{n+2p+\alp}{2n}}\\
&\leq C \left(\int_{\Sigma_\lambda^-}\left((C_0-\zeta)_+^{\gamma-1}(y)\right)^{\frac{n}{2p+\alp}}dy\right)^{\frac{2p+\alp}{n}}\left(\int_{\Sigma_\lambda^-}|w_\lambda(y)|^{\frac{2n}{n-2p-\alp}}dy\right)^{\frac{n-2p-\alp}{2n}}.
\end{split}
\end{align}
Recall that $\Sigma_\lambda^-\subset ((\Sigma_\lambda\setminus \Sigma_{\lambda_0 -\delta})\cap B_R)\cup (\Sigma_\lambda \setminus  B_R)$ and $-\zeta$ is bounded above, we can choose $\delta$ sufficiently small such that
$$ C \left(\int_{{\Sigma_\lambda^-}\cap B_R}\left((C_0-\zeta)_+^{\gamma-1}(y)\right)^{\frac{n}{2p}}dy\right)^{\frac{2p}{n}}\leq \frac{1}{2}. $$
From this and (\ref{p1}), we have
$$ C \left(\int_{\Sigma_\lambda^-}\left((C_0-\zeta)_+^{\gamma-1}(y)\right)^{\frac{n}{2p}}dy\right)^{\frac{2p}{n}}\leq \frac{1}{2}. $$
Now (\ref{x21}) implies that $\|w_\lambda\|_{L^{\frac{2n}{n-2p-\alp}}(\Sigma_\lambda^-)}=0$ and therefore $\Sigma_\lambda^-$ must be measure zero.

This verifies (\ref{x11}). Thus we must have $w_{\lambda_0}\equiv0$.

\par $\mathbf{Step~ 3}$: We show that $\frac{\partial \zeta}{\partial x_1}<0$ for $x\in\Sigma_{\lambda_0}$.

In fact, from the definition of $\lambda_0$ we have for any $\lambda<\lambda_0$,
\begin{equation}\label{z1}
w_\lambda(x)> 0, ~~~~~x\in \Sigma_\lambda.
\end{equation}
Simple calculation gives that for any $x\in T_\lambda$ with $\lambda<\lambda_0$,
\begin{equation*}
\begin{split}
\zeta_{x_1}(x)=&R_{2p+\alp,n}(n-2p-\alp)\int_{\mathbb{R}^n}\frac{(C_0-\zeta)_+^\gamma(y)(x_1-y_1)}{|x-y|^{n-2p+1-\alp}}dy\\
=&R_{2p+\alp,n}(n-2p-\alp)\int_{\Sigma_\lambda}\frac{\left((C_0-\zeta)_+^\gamma(y)-(C_0-\zeta_\lambda)_+^\gamma(y)\right)(x_1-y_1)}{|x-y|^{n-2p+1}}dy\\
<&0,
\end{split}
\end{equation*}
where the last inequality follows from (\ref{z1}). Thus the claim holds.

Since the problem is invariant with respect to rotation, we can take any direction as the $x_1$ direction. Hence we have that $\zeta$ is radially symmetric about some $x_0\in\mathbb{R}^n$ and $\frac{\partial \zeta}{\partial r}>0$ where $r=|x-x_0|$.
$\hfill\square$

Actually, we may also prove Theorem \ref{T3}  by applying  moving plane method to the function $(-\Delta)^{p-1+\frac{\alpha}{2}}u$, after asymptotic behaviors at infinity of this function and its first-order derivatives are established.

\end{document}